\documentclass[12pt]{amsart}

\usepackage{amssymb}

\newcommand{\ol}{\overline}

\newcommand{\al}{\alpha}

\newcommand{\dt}{\delta}

\newcommand{\gm}{\gamma}

\newcommand{\Dt}{\Delta}

\newcommand{\Tht}{\Theta}

\newcommand{\diag}{\mathop{\rm diag}}
\newcommand{\rank}{\mathbb{\rm rank}}
\newcommand{\ord}{\mathbb{\rm ord}}
\newcommand{\Cof}{\mathop{\rm Cof}}
\newcommand{\cof}{\mathop{\rm cof}}
\newcommand{\const}{\mathop{\rm Const}}

\newcommand{\tm}{\times}

\newcommand{\sbs}{\subset}

\newcommand{\wdt}{\widetilde}

\newcommand{\iy}{\infty}

\newcommand{\bR}{\mathbb{R}}

\newcommand{\bT}{\mathbb{T}}
\newcommand{\bC}{\mathbb{C}}
\newcommand{\bQ}{\mathbb{Q}}

\newcommand{\bF}{\mathbb{F}}

\newcommand{\calM}{\mathcal{M}}
\newcommand{\calU}{\mathcal{U}}
\newcommand{\calP}{\mathcal{P}}
\newcommand{\calA}{\mathcal{A}}
\newcommand{\calW}{\mathcal{W}}

\newcommand{\mbA}{\mathbf{A}}

\newcommand{\mbu}{\mathbf{u}}
\newcommand{\mbv}{\mathbf{v}}
\newcommand{\mbU}{\mathbf{U}}
\newcommand{\mbV}{\mathbf{V}}
\newcommand{\mbo}{\mathbf{0}}
\newcommand{\mbl}{\mathbf{1}}

\newcommand{\mbB}{\mathbf{B}}

\newtheorem{theorem}{\bf  Theorem}
\newtheorem{lemma}{\bf  Lemma}

\newtheorem{corollary}{\bf \sc Corollary}
\newtheorem{remark}{ \sc Remark}

\begin{document}


\begin{center}
On compact wavelet matrices of rank $m$ and of order and degree $N$
\\[5mm]

 L.~Ephremidze  and  E.~Lagvilava
        \end{center}

\vskip+0.5cm

 \noindent {\small {\bf Abstract.} A new parametrization (one-to-one onto map) of
  compact wavelet matrices of rank $m$ and of order and degree $N$ is proposed in terms of coordinates
  in the Euclidian space $\bC^{(m-1)N}$. The developed method depends  on
  Wiener-Hopf factorization of corresponding unitary matrix functions and allows to construct
  compact wavelet matrices efficiently. Some applications of the proposed method are discussed.}

 \vskip+0.2cm\noindent  {\small {\bf Keywords:} Wavelet matrices,
paraunitary matrix polynomials, Wiener-Hopf factorization.}

\vskip+0.2cm \noindent  {\small {\bf  AMS subject classification
(2010):} 42C40}

\vskip+0.5cm

\section{Introduction}

Let $\bC$ be the set of complex numbers, and let $\bF\sbs\bC$ be a
subfield invariant under the complex conjugation. A formally
infinite matrix $\calA$ with $m$ rows
\begin{equation}
\label{WM1} {\mathcal A}= \left(\begin{matrix} \cdots&
a^0_{-1}&a^0_{0}&a^0_{1}&a^0_{2}&\cdots\\[1mm]
\cdots& a^1_{-1}&a^1_{0}&a^1_{1}&a^1_{2}&\cdots\\
&\vdots&\vdots\\
\cdots& a^{m-1}_{-1}&a^{m-1}_{0}&a^{m-1}_{1}&a^{m-1}_{2}&\cdots\\
\end{matrix}\right),\;\;\;a^i_j\in\bF,
\end{equation}
is called a {\em wavelet matrix} (of rank $m$) if its rows satisfy
the so called {\em shifted orthogonality condition}:
\begin{equation}
\label{Qc2}
 \sum_{k=-\iy}^\iy a^i_{k+mj}\,\ol{a^r_{k+ms}}=\dt_{ir}\dt_{js}
\end{equation}
($\dt$ stands for the  Kronecker delta). Such matrices are a
generalization of ordinary $m\tm m$ unitary matrices and they play
the crucial role in the theory of wavelets \cite{RW} and multirate
filter banks \cite{Vai}. The reason for us to work with $\bF$
instead of just $\bC$ is that it will allow more flexibility in
 applications, since the proposed proofs and discussions apply
to a whole range of fields including the set of real numbers,
rational numbers, as well as algebraic extensions of rational numbers. A reader not
concerned with general fields may  assume that $\bF=\bC$.

In the {\em polyphase representation} of matrix $\calA$,
\begin{equation}
\label{PR3} {\mathbf A}(z)=\sum_{k=-\iy}^\iy A_kz^k\,,
\end{equation}
where ${\mathcal A}=(\cdots, A_{-1},A_0,A_1,A_2,\cdots)$ is the
partition of $\calA$ into block $m\tm m$ matrices
$A_k=(a^i_{km+j})$, $0\leq i,j\leq m-1$, the condition (\ref{Qc2})
is equivalent to
\begin{equation}
\label{Pp35} {\mathbf A}(z)\wdt{\mbA}(z)=I_m,
\end{equation}
where $\wdt{\mbA}(z)=\sum_{k=-\iy}^\iy A_k^*z^{-k}$ is the {\em
adjoint} to $\mbA(z)$ ($A^*:=\ol{A}^T$ is the Hermitian conjugate
and $I_m$ stands for the $m\tm m$ unit matrix). This is easy to
see as (\ref{Qc2}) can be written in the block matrix form
$\sum_{k=-\iy}^\iy A_k A^*_{l-k}=\dt_{l0}I_m$.

Our notion of a wavelet matrix is weaker than usual. So, as the
orthogonal basis of $L^2(\bR)$ can be developed (see \cite[Ch-s 4,
5]{RW}), also the {\em linear condition}
$\mbA(1)\,\mbl=\sqrt{m}\,e_1$, where $\mbl=(1,1,\ldots,1)^T$ and
$e_1=(1,0,\ldots,0)$, must be satisfied. In our consideration, the
linear condition is irrelevant. Instead, we require the  condition
\begin{equation}
\label{id4} {\mathbf A}(1)=I_m
\end{equation}
which much simplifies  the whole exposition. For any wavelet matrix
$\calA$, $\mbA(1)=\sum_{k=-\iy}^\iy A_k$ is unitary (see
(\ref{Pp35})), so that the conditions (\ref{Qc2}) and (\ref{id4})
will be satisfied for $\calA_0=\big(\mbA(1)\big)^{-1}\calA$. Note
also that for any wavelet matrix $\calA_0$ satisfying (\ref{id4})
and any unitary matrix $U$, the matrix $\calA=U\calA_0$ is also a
wavelet matrix satisfying ${\mathbf A}(1)=U$. Thus the additional
restriction (\ref{id4}) does not lose generality in the
description of wavelet matrices. Observe that the polyphase
representation of $U\calA$ is $\sum_{k=-\iy}^\iy UA_kz^k$
(obviously, there is one-to-one correspondence between the matrices
(\ref{WM1}) and their polyphase representations (\ref{PR3}) and
they are naturally identified).

We consider the compact wavelet matrices, which means that only a finite
number of coefficients in (\ref{WM1}) are non-zero (the corresponding
wavelet functions in $L^2(\bR)$ have then a compact support, and the
corresponding filters in signal processing applications are
physically realizable). Namely, compact wavelet matrices with
polyphase representation
\begin{equation}
\label{Mp5} {\mathbf A}(z)=\sum_{k=0}^N A_kz^k\,,
\end{equation}
where $A_N\not=0$, are called of ({\em rank} $m$ and) {\em order}
$N:=\ord(\calA)$ (in some works they are called of {\em genus}
$N+1$); we write $\calA\in\mathcal{W}\mathcal{M}(m,N,\bF)$. The
property (\ref{Pp35}) for the matrix polynomial (\ref{Mp5}) means that
$\mbA(z)$ is a {\em paraunitary} matrix function. Note that, in this
case, $\mbA(z)$ is usual unitary matrix for each
$z\in\bT:=\{z\in\bC:|z|=1\}$.

As the determinant of $\mbA(z)\in\mathcal{W}\mathcal{M}(m,N,\bF)$
is always monomial, i.e., $\det\mbA(z)=c\,z^d$ (since
$\det\mbA(z)\det\wdt{\mbA}(z)=1$), the positive integer $d$ is
defined to be the {\em degree} of $\calA$, $d=\deg(\calA)$. The
set of wavelet matrices of rank $m$, order $N$ and degree $d$ will
be denoted by $\mathcal{W}\mathcal{M}(m,N,d,\bF)$. It can be
proved that $\deg(\calA)\geq\ord(\calA)$ for any
$\calA\in\mathcal{W}\mathcal{M}(m,N,\bF)$  and
$\deg(\calA)=\ord(\calA)$ if and only if $\rank(A_0)=m-1$ (see
Lemma 1 below). As $A_0A_N^*=\mbo$ for each
$\calA\in\mathcal{W}\mathcal{M}(m,N,\bF)$ (since it is the $N$th
matrix coefficient of the product in (\ref{Pp35}){}) and
$A_N\not=\mbo$, we have $\rank(A_0)<m$. Thus $\deg(\calA)=N$
whenever $A_0$ has a maximal possible rank and it is assumed to be
the nonsingular case. Throughout the paper, we will intensively
study $\mathcal{W}\mathcal{M}(m,N,N,\bF)$.

Obviously, if we shift rows of  $\calA\in
\mathcal{W}\mathcal{M}(m,N,N,\bF)$ by arbitrary multiples of $m$,
then the  obtained $\calA_1$ is a wavelet matrix as well. However, if we
make these movements chaotically, then in general $\deg(\calA_1)$
does not coincide with $\ord(\calA_1)$  and $\calA_1$ becomes
singular.

 As $\calA$ and $U\calA$, where $U$ is nonsingular,
have the same rank, order and degree, we will additionally assume that
 $\calA\in\mathcal{W}\mathcal{M}(m,N,N,\bF)$ satisfies
(\ref{id4}) and the subset of such compact wavelet matrices will
be denoted by $\mathcal{W}\mathcal{M}_0(m,N,N,\bF)$. As has been
mentioned above, the linear condition (\ref{id4}) does not lose
generality in the description of
$\mathcal{W}\mathcal{M}(m,N,N,\bF)$.

Matrix polynomials $\mbV(z)\in\mathcal{W}\mathcal{M}_0(m,1,1,\bF)$
are called {\em primitive} wavelet matrices and they always have the
form $\mbV(z)=Q+Pz$, where $P$ and $Q$ are {\em complementary
$($orthogonal$)$ projections} on $\bF^m$ and
$\rank(P)=1\Leftrightarrow\rank(Q)=m-1$ (see \cite[Th. 3.1]{KT} or
Lemma 2).

Every $\mbA(z)\in\mathcal{W}\mathcal{M}_0(m,N,N,\bF)$ can be
uniquely factorized as (see \cite[Theorem 4.4.15]{RW} or Theorem
2)
\begin{equation}
\label{pf6} {\mathbf A}(z)=\prod_{j=1}^N \mbV_j(z)=\prod_{j=1}^N
\big(Q_j+P_jz\big),
\;\;\;\mbV_j(z)\in\mathcal{W}\mathcal{M}_0(m,1,1,\bF),
\end{equation}
where no consecutive operators $P_j$ are
 orthogonal to each other, $P_jP_{j+1}\not=0$. This factorization
gives  rise to the map (parametrization)
 \begin{equation}
\label{par1} \underbrace{\bF \mathrm{P}^{m-1}\tm\bF
\mathrm{P}^{m-1}\tm\ldots\tm\bF \mathrm{P}^{m-1}}_N\supset
\mathcal{B} \longleftrightarrow
\mathcal{W}\mathcal{M}_0(m,N,N,\bF),
\end{equation}
where $\bF\mathrm{P}^{m-1}$ is the projective space
$\bF^m/(\bF\backslash \{0\})$ (the space of one-dimensional
subspaces of $\bF^m$) which is one-to-one and onto but defined
only on the subset $\mathcal{B}$ where consecutive directions are
not orthogonal, i.e. some singular points are excluded from the
set of parameters. Until now it has been
 the only known simple way of constructing $\calA\in\mathcal{W}\mathcal{M}_0(m,N,N,\bF)$
 for arbitrary $N$: choosing nonzero column
 vectors $v_j\in\bF^m$,  $j=1,2,\ldots,N$, such
 that $v_j\not\perp v_{j+1}$, taking the corresponding projections
  $P_j=v_j(v_j^*v_j)^{-1}v_j^*$ and primitive wavelet matrices
 $\mbV_j=I_m-P_j+P_jz$, and constructing the product (\ref{pf6})
that belongs to $\mathcal{W}\mathcal{M}_0(m,N,N,\bF)$.

\begin{remark}
The representation $(\ref{pf6})$ resembles the factorization for
ordinary polynomials $p(z)=\prod_{k=1}^N(z-z_k)$ by using their
roots. However the following remarkable fact should be emphasized:
to find the roots of $p(z)$ we have to go, in general, to a larger
field than the field of coefficients of $p(z)$. For the
factorization in $(\ref{pf6})$, the coefficients of each factor
$\mbV_j$ belong to the same field $\bF$.
\end{remark}

Now we turn to our contribution in the study of compact wavelet
matrices. This approach has been developed during a search for a
new matrix spectral factorization algorithm \cite{JLE}, and it
makes possible to parameterize $\calW\calM(m,N,N,\bF)$ directly in
terms of points in $\bF^{(m-1)N}$. Moreover, compact wavelet
matrices can be constructed in a more efficient way than by taking
the products (\ref{pf6}). Other advantages and possible
applications of the proposed method are mentioned in the course of discussions
below.

Recall that we identify $\calA\in\calW\calM(m,N,N,\bF)$ with its
polyphase representation (\ref{Mp5}). The linear condition we have
introduced is (\ref{id4}) and such a subclass is denoted by
$\calW\calM_0(m,N,N,\bF)$. We further require that the last row of
$A_N$ be not the zero vector and denote such a subclass by
$\calW\calM_1(m,N,N,\bF)$. This is done without loss of generality
as $A_N\not=0$ and we can interchange the rows of $\calA$, if
necessary. In particular, for any $\calA\in\calW\calM(m,N,N,\bF)$
there exist constant unitary $m\tm m$ matrices $U_0$ and $U_1$
such that $U_0\calA U_1\in\calW\calM_1(m,N,N,\bF)$ (by  right
multiplication  we assume that $\calA U_1=(A_0U_1,
A_1U_1,\ldots,A_NU_1)$, i.e. the polyphase matrix of $\calA U_1$
is $\mbA(z)U_1$). Indeed, we can interchange the rows of $\calA$,
if necessary, by multiplication by $U_0$ and then take
$U_1=\big(U_0\mbA(1)\big)^{-1}$. Hence, in what follows, we
parameterize $\calW\calM_1(m,N,N,\bF)$.

Let   $\calP_N^+[\bF]:=\big\{\!\sum_{k=0}^N c_k
z^k:c_0,c_1,\ldots,c_N$ $\in\bF\big\}$ be the set of polynomials
with coefficients from $\bF$ and
$\calP_N^-[\bF]:=\big\{\sum_{k=1}^N c_k
z^{-k}:c_1,c_2,\ldots,c_N\in\bF\big\}$ which can be naturally
identified with $\bF^{N+1}$ and $\bF^N$, respectively (note that
$\calP_N^+[\bF]\cap\calP_N^-[\bF]=\{0\}$ according to our
notation). Sometimes $[\bF]$ is omitted because it is clear
from the context.

 Let $\calP\calU_1(m,N,\bF)$ be a set of $m\tm m$ paraunitary matrix
functions $U(z)$,
\begin{equation}
\label{UPp1} \wdt{U}(z)U(z)=I_m,
\end{equation}
of the special form
\begin{equation}
\label{U6}
U(z)=\begin{pmatrix}u_{11}(z)&u_{12}(z)&\cdots&u_{1m}(z)\\
                 u_{21}(z)&u_{22}(z)&\cdots&u_{2m}(z)\\
           \vdots&\vdots&\vdots&\vdots\\
           u_{m-1,1}(z)&u_{m-1,2}(z)&\cdots&u_{m-1,m}(z)\\[3mm]
           \wdt{u_{m1}}(z)&\wdt{u_{m2}}(z)&\cdots&\wdt{u_{mm}}(z)\\
           \end{pmatrix},\;\;\;u_{ij}(z)\in\calP_N^+[\bF],
\end{equation}
$1\leq i,j\leq m$, with determinant 1,
\begin{equation}
\label{Udet3} \det U(z)=1,
\end{equation}
and with the linear restriction
\begin{equation}
\label{UIm4} U(1)=I_m,
\end{equation}
and such that not all polynomials $u_{mj}$ (these are the adjoint
polynomials of the entries in the last row) are zero at the origin,
\begin{equation}
\label{Ulr5} \sum_{j=1}^m|u_{mj}(0)|>0.
\end{equation}
Then for each $\mbA(z)\in \calW\calM_1(m,N,N,\bF)$, we have
$U(z)=\diag[1,\ldots,1,z^{-N}]\mbA(z)\in \calP\calU_1(m,N,\bF)$
(the last row is multiplied by $z^{-N}$), and conversely. Thus
there is a simple one-to-one correspondence
\begin{equation}
\label{W1U1} \calW\calM_1(m,N,N,\bF)\longleftrightarrow
\calP\calU_1(m,N,\bF),
\end{equation}
and we parameterize the latter class using the following

\begin{theorem}
Let $N\geq 1$. For any  Laurent matrix polynomial  $F(z)$ of the
form
\begin{equation}
\label{F93}
F(z)=\begin{pmatrix}1&0&0&\cdots&0&0\\
          0&1&0&\cdots&0&0\\
           0&0&1&\cdots&0&0\\
           \vdots&\vdots&\vdots&\vdots&\vdots&\vdots\\
           0&0&0&\cdots&1&0\\
           \zeta_{1}(z)&\zeta_{2}(z)&\zeta_{3}(z)&\cdots&\zeta_{m-1}(z)&1
           \end{pmatrix},\;\;\;\zeta_j(z)\in\calP_N^-[\bF],
\end{equation}
there exists a unique
\begin{equation}
\label{U94} U(z)\in\mathcal{P}\mathcal{U}_1(m,N,\bF)
\end{equation}
such that
\begin{equation}
\label{Pr95} F(z)U(z)\in\calP^+(m\tm m).
\end{equation}

Conversely, for each $U(z)$ satisfying $(\ref{U94})$, there exists
a unique matrix function $F(z)$ of the form $(\ref{F93})$ such
that $(\ref{Pr95})$ holds.
\end{theorem}

$\calP^+(m\tm m)$ stands for the set of $m\tm m$ matrix
polynomials with entries from $\calP^+$.

\begin{remark}
 The relation $(\ref{Pr95})$ written in the equivalent form
 $U(z)=F^{-1}(z)M(z)$, where $M(z)\in \calP^+(m\tm m)$ and
  each $\zeta_j$ is replaced by $-\zeta_j$ in $(\ref{F93})$ to get
 $F^{-1}(z)$, means that we have the Wiener-Hopf factorization of
 the unitary matrix function $U(t)=U(z)|_{z=t}$ defined on $\bT$,
 since $M(z)$ and $F^{-1}(z)$ are analytic regular matrix functions inside
 and outside $($including infinity$)$ $\bT$, respectively.
 Although we heavily exploited this idea in our previous works on the
 factorization of matrix functions \cite{Proc1998}, \cite{JLE},
 from which this paper has stemmed out,  we tried
 to reduce to a minimum the application of Wiener-Hopf
 factorization theory in the presented proofs and discussions.
 This allowed us to transfer the obtained results to arbitrary
 subfields $\bF\sbs\bC$ which are invariant  under conjugate
 operation, and thus extend the area of their applications.
\end{remark}

The proof of Theorem 1 makes it possible  to explicitly construct
the corresponding $U(z)$  for a given
$F(z)$,   and vice versa.
  Namely, let the functions in (\ref{F93}) be
\begin{equation}
\label{zeta1} \zeta_i(z)=\sum_{k=1}^N \gm_{ik}z^{-k},\;\;\;
i=1,2,\ldots,m-1,
\end{equation}
and let $\Tht_i$ be the $(N+1)\tm(N+1)$ Hankel-like matrix
\begin{equation}
\label{tht1}
\Tht_i=\begin{pmatrix} 0& \gm_{i1}& \gm_{i2}&\cdots& \gm_{i,N-1}& \gm_{iN}\\
         \gm_{i1}& \gm_{i2}& \gm_{i3}&\cdots& \gm_{iN}&0\\
         \gm_{i2}& \gm_{i3}& \gm_{i4}&\cdots&0&0\\
        \cdot&\cdot&\cdot&\cdots&\cdot&\cdot\\
         \gm_{iN}&0&0&\cdots&0&0\end{pmatrix},
         \;\;\;i=1,2,\ldots,m-1,
\end{equation}
and
\begin{equation}
\label{Dt1} \Dt=\sum_{i=1}^{m-1}\Tht_i\ol{\Tht_i}+I_{N+1}.
\end{equation}
 Assume also that $B_i$ is the first column of $\Tht_i$,
$i=1,2,\ldots,m-1$, and $B_m=(1,0,\ldots,0)^T$. Let
$X_j=(\al_{j0},\al_{j1},\ldots,\al_{jN})^T\in\bF^{N+1}$ be the
solution of the following linear system of algebraic equations
\begin{equation}
\label{517} \Dt\cdot X=B_j,\;\;j=1,2,\ldots,m
\end{equation}
($\Dt$ is positive definite and has a {\em displacement structure
of rank $m$}, therefore O$(mN^2)$ operations are required for its
solution instead of traditional O$(N^3)$; see
\cite[Appendix]{JLE}),  and let
\begin{gather}
v_{mj}(z)=\sum_{k=0}^N \al_{jk}z^{-k},\;\;\;j=1,2,\ldots,m,\\
v_{ij}(z)=\big[\wdt{\zeta_i}(z)v_{mj}(z)\big]^+-\dt_{ij},\;\;\;
i=1,2,\ldots,m-1,
\end{gather}
where $[\;\cdot\;]^+$ stands for the projection operator,
$\left[\sum_{k=-N}^N c_k z^k\right]^+= \sum_{k=0}^N c_k z^k$. Then
\begin{equation}
\label{V-1} U(z)=V(z)\big(V(1)\big)^{-1},
\end{equation}
where $V(z)$ is the $m\tm m$ Laurent polynomial matrix $
V(z)=\left(v_{ij}(z)\right)_{i,j=1}^m$, will be the desired matrix
polynomial (\ref{U94}) as  proved in Sect. 4.

For a given $U(z)$, the corresponding $F(z)$ can be also
explicitly computed as follows (see Sect. 4 for the proof): if
$u_{mj}(0)\not=0$, then
\begin{equation}
\label{zeta24}
\zeta_i(z)=\left[\frac{\wdt{u}_{ij}(z)}{u_{mj}(z)}\right]^-,
\;\;\;\;i=1,2,\ldots,m-1,
\end{equation}
where $[\;\cdot\;]^-$ stands for the projection operator:
$\left[\sum_{k=-N}^\iy c_k z^k\right]^-= \sum_{k=-N}^{-1} c_k
z^k$, and under $\frac 1{u_{mj}(z)}$ its formal series expansion
in a neighborhood of $0$ is assumed. Note that we need to take
only the first $N+1$ coefficients in this expansion in order to
compute $\zeta_i$s.

In consequence, as the set of matrix polynomials of the form
(\ref{F93}) can be easily identified with $\bF^{(m-1)N}$, we have
the following diagram of one-to-one and onto maps
\begin{equation}
\label{PAR} \calW\calM_1(m,N,N,\bF)\longleftrightarrow
\calP\calU_1(m,N,\bF)\longleftrightarrow \underbrace{\bF^N\tm\bF^N
\tm\ldots\tm\bF^N}_{m-1}\,,
\end{equation}
i.e. a complete parametrization of nonsingular compact wavelet
matrices, which can be effectively realized.

The paper is organized as follows. In the next section  we provide
all notation and definitions used throughout the paper. Most of
them have already been introduced in the current section but we
collect them together for the convenience of reference. In Sect. 3 we
reprove  the factorization (\ref{pf6}) and related facts as it
seems that our approach is somewhat  simpler than the traditional one.
The main result of our paper, Theorem 1, is proved in Sect. 4. The
essential part of this theorem, the existence of the paraunitary matrix
function (\ref{U94}) (without showing the property (\ref{Ulr5})),
has already been published (see \cite[Th. 1]{JLE}) for the case
$\bF=\bC$.\footnote{This constructive proof is a core of a new
matrix spectral factorization method published in \cite{JLE}. This
method is currently patent pending.} It is sufficient to note that
the proof goes through without any change if we replace $\bC$ by
$\bF$. However, for readers' convenience, we  included
the simplified version of this proof in the present paper. (A
minor inaccuracy of the discussion in \cite{JLE} between the formulas
(50) and (51) has also been corrected.)

In the remaining two sections we consider some possible
applications of our method. In particular, in Sect 5, we solve
the following important problem: given the first row of a wavelet
matrix $\calA$ (sometimes called the {\em scaling vector} or {\em
low-pass filter}), how to find the remaining rows of $\calA$
(called the {\em wavelet vectors} or {\em high-pass filters}). The
solution of this problem is well known and there exists an
appropriate algorithm for construction of wavelet vectors (see
e.g. \cite[p. 64]{RW}). However this algorithm requires too many
matrix multiplications which might cause the round-off problems
for  large $m$ and $N$. We propose another algorithm for solution
of this problem which uses our approach. A complete comparison
of both algorithms based on numerical simulations will be
provided soon.

Sect. 6 contains some discussions about other possible
applications of the proposed parametrization of compact wavelet
matrices.

\section{Notation and Preliminaries}

The fields of rational, real, and complex numbers are denoted by
$\bQ, \bR$, and $\bC$, respectively, and $\bF$ stands for a
subfield of $\bC$ which is invariant under the complex
conjugation, $a\in\bF\Rightarrow \ol{a}\in\bF$.

$\bT:=\{z\in\bC:|z|=1\}$, $\bT_+:=\{z\in\bC:|z|<1\}$, and
$\bT_-:=\{z\in\bC:|z|>1\}\cup\{\infty\}$.

Let ${\bf e}_1,{\bf e}_2,\ldots,{\bf e}_m$ be the standard basis
of $\bF^m$, e.g., ${\bf e}_1=(1,0,\ldots,0)$. The usual scalar
product and norm in the space $\bF^m$ are denoted by
$\langle\,\cdot,\cdot\rangle_m$ and $\|\cdot\|_m$, respectively
(they are independent of $\bF$ as it is a subfield of $\bF$).

Let $\bF^{m\tm n}$ be the set of $m\tm n$ matrices with entries
from a field $\bF$. For $M\in \bF^{m\tm n}$, let $M^*=\ol{M}^T\in
\bF^{n\tm m}$ be the Hermitian conjugate of $M$.
$\diag[a_1,a_2,\ldots,a_m]$ is an $m\tm m$ diagonal matrix with
corresponding entries on the diagonal and
$I_m=\diag[1,1,\ldots,1]$ is the $m\tm m$ unit matrix. $U\in
\bF^{m\tm m}$ is called unitary, $U\in\calU_m[\bF]$, if
$UU^*=U^*U=I_m$.

If $M\in \bF^{m\tm m}$ and $a$ is an entry of $M$, then
$\cof(a)\in\bF$ is the cofactor of $a$ and $\Cof(M)\in\bF^{m\tm
m}$ is the cofactor matrix, so that $M^{-1}=\frac 1{\det
M}\big(\Cof(M)\big)^T$ if $M$ is nonsingular. The same notation
will be used for matrix functions.

A matrix $P\in \bF^{m\tm m}$ (considered as a linear map from
$\bF^m$ to $\bF^m$) is called the (orthogonal) projection if it is
{\em self-adjoint}, $P=P^*$, and {\em idempotent}, $P^2=P$,
however we will reduce this definition to a single formula:
$P=PP^*$. Note that if $P\in \bF^{m\tm m}$ is a projection of rank
$r\leq m$, then there exists $U\in\bC^{m\tm r}$ with orthonormal
columns, $U^*U=I_r$, such that
\begin{equation}
\label{2Pu} P=UU^*
 \end{equation}
(indeed, the factorization (\ref{2Pu}) with $U\in\bC^{m\tm r}$ of
full rank $r$ exists since $P$ is non-negative definite, and
$UU^*UU^*=UU^*$ $\Rightarrow$ $U^*UU^*U U^*U=U^*UU^*U$
$\Rightarrow$
$(U^*U)^{-1}\tm(U^*U)(U^*U)(U^*U)(U^*U)^{-1}=(U^*U)^{-1}
(U^*U)(U^*U)(U^*U)^{-1}$ $\Rightarrow$ $U^*U=I_r$).

If $P$ and $Q$ from $\bF^{m\tm m}$ are projections and $P+Q=I_m$,
then they are called complementary to each other. Obviously,
$PQ^*=Q^*P={\mbo}$ for complementary projections (as
$P(I_m-P)^*=P-PP^*=\mbo$), and if $P$ is a projection, then
$I_m-P$ is also the projection complementary to $P$.

$\calP[\bF]$ denotes the set of Laurent polynomials with
coefficients from a field $\bF$, and
$\calP_N[\bF]:=\{\sum_{k=-N}^N c_kz^k:c_k\in \bF,
k=-N,\ldots,N\}$. If we write just $\calP$, the field of
coefficients will be clear from the context. $\calP^+\sbs\calP$ is
the set of polynomials (with non-negative powers of $z$,
$\sum_{k=0}^N c_kz^k\in \calP^+$) and $\calP^-\sbs\calP$ is the
set of Laurent polynomials with negative powers of $z$,
$\sum_{k=1}^N c_kz^{-k}\in \calP^-$. We emphasize that, according to our notation,
 constant functions belong only to $\calP^+$  so
that $\calP^+\cap \calP^-=\{0\}$. Let also
$\calP^\pm_N=\calP^\pm\cap\calP_N$.

$\calP(m\tm n)$ denotes the set of $m\tm n$ (polynomial) matrices
with entries from $\calP$, and  the sets $\calP^+(m\tm n)$,
$\calP^-_N[\bF](m\tm n)$, etc. are defined similarly.  The
elements of these sets, $P(z)=\big(p_{ij}(z)\big)$, are called
(polynomial) matrix functions. When $n=1$ we have the vector
functions and such classes are denoted by $\calP(m)$ instead of
$\calP(m\tm 1)$. When we speak about the continuous maps between these
sets, we mean that they are equipped with a usual topology.

For $p(z)=\sum_{k=-N}^N c_kz^k \in\calP$, let
$\wdt{p}(z)=\sum_{k=-N}^N \ol{c}_kz^{-k}$ and for
$P(z)=[p_{ij}(z)]\in \calP(m\tm n)$ let
$\wdt{P}(z)=[\wdt{p}_{ij}(z)]^T\in\calP(n\tm m)$.  Note that
$\wdt{P}(z)=(P(z))^*$ when $z\in\bT$. Thus usual relations for
adjoint matrices like
$\wdt{P_1+P_2}(z)=\wdt{P}_1(z)+\wdt{P}_2(z)$,
$\wdt{P_1P_2}(z)=\wdt{P}_2(z)\wdt{P}_1(z)$, etc. hold.

We employ also the additional notation of the sets
$\wdt{\calP_N^+}[\bF]:=\{\wdt{p}(z):p(z)\in\calP_N^+[\bF]\}$,
which is an extension of $\calP_N^-[\bF]$, and
$\calP_N^\oplus[\bF](m):=
\underbrace{{\calP_N^+}[\bF]\tm\ldots\tm{\calP_N^+}[\bF]\tm
\wdt{\calP_N^+}[\bF]}_m$.

A polynomial matrix (\ref{Mp5}), where $A_k\in\bF^{m\tm m}$,
$k=1,2,\ldots,N$, is called paraunitary, we write
$\mbA(z)\in\calP\calU(m,N,\bF)$, if (\ref{Pp35}) holds, which is
equivalent to $\sum_{k=0}^{N-l}A_kA^*_{k+l}=\dt_{l0}I_m$,
$l=0,1,\ldots,N$. This means that ${\mathcal
A}=(A_0,A_1,A_2,\cdots,A_N)\in\bF^{m\tm m(N+1)}$ is a wavelet
matrix and we do not make distinction between $\calA$ and its its
polyphase representation (\ref{Mp5}). If $A_N\not=\mbo$ in the
representation (\ref{Mp5}), then we say that the wavelet matrix
$\calA\equiv \mbA(z)$ has rank $m$ and order $N$, we write
$\mbA(z)\in\calW\calM(m,N,\bF)$. The degree of a wavelet matrix
$\mbA(z)$ is the order of a monomial $\det\mbA(z)$, i.e.
$\det\mbA(z)=c\cdot z^{\deg(\calA)}$. The set of wavelet matrices
of rank $m$, order $N$, and degree $d$ will be denoted by
$\mathcal{W}\mathcal{M}(m,N,d,\bF)$. In addition,
$\mathcal{W}\mathcal{M}_0(m,N,d,\bF)$ is the subset of wavelet
matrices which satisfy (\ref{id4}), and
$\mbA(z)\in\mathcal{W}\mathcal{M}_1(m,N,d,\bF)\sbs
\mathcal{W}\mathcal{M}_0(m,N,d,\bF)$ if the last row of $A_N$
differs from the zero row vector.

Recall also that $\calP\calU_1(m,N,\bF)$ denotes the set of
specific paraunitary matrix functions which satisfy
(\ref{UPp1})--(\ref{Ulr5}).

For power series $f(z)=\sum_{k=-\iy}^\iy c_k z^k$ and $N\geq 1$,
let $[f(z)]^-$, $[f(z)]^+$, $[f(z)]^-_N$, and $[f(z)]^+_N$,
denote respectively $\sum_{k=-\iy}^{-1} c_k z^k$,
$\sum_{k=0}^\iy c_k z^k$, $\sum_{k=-N}^{-1} c_k z^k$, and
$\sum_{k=0}^N c_k z^k$ and the corresponding functions are assumed
under these expressions if the convergence domains of these power
series are clear from the context (the dependence of coefficients
$c_k$ on the function $f$ is sometimes expressed by
$c_k\{f\}=c_k${}). Obviously $[f]^-=f-[f]^+$
 under these notation.

A matrix function $S(z)\in\calP[\bC](m\tm m)$ is called positive
definite if $S(z)$ is such ($XS(z)X^*>0$ for each
$0\not=X\in\bC^{1\tm m}$) for almost every $z\in\bT$. The
polynomial matrix spectral factorization theorem (see e.g.
\cite{EJL2}) asserts that every positive definite
$S(z)\in\calP_N(m\tm m)$ can be factorized as
\begin{equation}
 S(z)=S^+(z)\wdt{S^+}(z),\;\;\;z\in\bC\backslash \{0\},
 \label{PSF}
 \end{equation}
where $S^+\in\calP^+_N(m\tm m)$ and $\det S^+(z)\not=0$ for each
$z\in \bT_+$ (consequently, $\wdt{S^+}(z)$ is analytic and
invertible in $\bT_-$), and the representation (\ref{PSF}) is
unique in a sense that if $ S(z)=S^+_1(z)\wdt{S^+_1}(z)$ is
another spectral factorization of $S(z)$, then  $S^+_1(z)=S^+(z)U$
for some unitary $U\in\calU_m$.

\section{Wavelet  Matrix Factorization Theorem}

The material of this section is mostly well known, however we
provide compact proofs of the given statements making emphasis on
arbitrariness   of a field $\bF$.

\begin{lemma} Let
\begin{equation}
\label{Mp52} {\mathbf A}(z)=\sum_{k=0}^N
A_kz^k\,,\;\;\;A_k\in\bF^{m\tm m},\;\;A_N\not=0,
\end{equation}
be a wavelet matrix of rank $m$, order $N$ and degree $d$,
$\mbA(z)\in\calW\calM(m,N,d,\bF)$. Then

a$)$\; $d\geq N$;

b$)$ $d=N$ if and only if $\rank(A_0)=m-1$.
\end{lemma}

\begin{proof}
We have
$$
\sum_{k=0}^N
 A^*_kz^{-k}= \widetilde{\mbA}(z)=\mbA^{-1}(z)=\frac1{\det \mbA(z)}\big(\Cof \mbA(z)\big)^T=
 cz^{-d} \big(\Cof \mbA(z)\big)^T,
$$
and since  $ \Cof\mbA(z)\in\calP^+(m\tm m)$ and $A_N\not=0$, the
equation $\sum\limits_{k=0}^N
 A^*_kz^{-k}\!= \frac{c}{z^{d}} \big(\Cof \mbA(z)\big)^T$ implies a). Using
 the same reasoning, we conclude that
 $\rank(A_0)<m-1\Longleftrightarrow \Cof(A_0)=\mbo\Longleftrightarrow
d>N $ (note that if
 $\Cof\mbA(z)=\sum_{k=0}^N C_kz^k$, then
 $C_0=\Cof(A_0)$).
 Hence b) follows as well (recall that $\rank(A_0)\not=m$ since
 $A_0A^*_N=0$).
\end{proof}

The following lemma will be used in the sequel only for $d=1$.

\begin{lemma}{\rm (cf. \cite[Th. 3.1]{KT}).}
Let $\mbV(z)=Q+Pz$, $P, Q\in\bF^{m\tm m}$, be a  matrix
polynomial. Then it is a wavelet matrix of rank $m$, order  1 and
degree $d$ satisfying $V(1)=I_m$, $\mbV(z)\in
\calW\calM_0(m,1,d,\bF)$, if and only if $P$ and $Q$ are
complementary projections on $\bF^m$ and $\rank(P)=d$.
\end{lemma}

\begin{proof}
The first part of the lemma can be proved by observation that
$(Q+Pz)(Q^*+P^*z)=I_m$ and $\mbV(1)=I_m$ $\Longleftrightarrow$
$PQ^*=0,\;PP^*+QQ^*=I_m$, and $P+Q=I_m$ $\Longleftrightarrow $
 $P$ and $Q$ are complementary projections.

Now the property of $P$ to be a projection is equivalent to the
existence of $U\in\bC^{m\tm d}$, $\rank(U)=\rank(P)=d$, with
orthonormal columns, $U^*U=I_d$, such that (\ref{2Pu}) holds. Let
$U_0\in\bC^{m\tm m}$ be any unitary matrix which completes the
columns of $U$. Then
$$
U_0^*V(z)U_0=U_0^*(I-UU^*+UU^*z)U_0=\diag[\underbrace{z,z,\ldots,z}_d,
\underbrace{1,1,\ldots,1}_{m-d}]
$$
whose determinant is $z^d$. Hence the lemma is proved.
\end{proof}

\begin{theorem}{\rm (see \cite[Th. 4.4.15]{RW}).}
Let $\mbA(z)$ be a wavelet matrix of rank $m$ and of order and
degree $N$ satisfying the linear condition $(\ref{id4})$,
$\mbA(z)\in\calW\calM_0(m,N,N,\bF)$. Then there exists a unique
factorization $(\ref{pf6})$.
\end{theorem}

\begin{proof}
By virtue of Lemma 1, $\rank(A_0)=m-1$. Let a nonzero column
vector $v\in\bF^m$ be orthogonal to the columns of $A_0$,
$v^*A_0=0$, and $\mbV_1(z)=I_m-P+Pz$, where $P=v(v^*v)^{-1}v^*$ is
the projection. Then $\mbV_1(z)\in\calW\calM_0(m,1,1,\bF)$ and
$\mbV_1(z)$ divides $\mbA(z)$ on the left (the quotient
$\mbB(z):=\mbV_1^{-1}(z)\mbA(z)=\wdt{\mbV_1}(z)\mbA(z)\in\calP^+(m\tm
m)$ as $P^*A_0=0$).

As $\deg(\mbB)=\deg(\mbA)-\deg(\mbV_1)=N-1$ and
$\deg(\mbB)\geq\ord(\mbB)$ by virtue of Lemma 1, we have
$\ord(\mbB)=N-1$ (obviously, it cannot be less than $N-1$ as
$\mbV_1(z)\mbB(z)=\mbA(z)${}). Hence
$\mbB(z)\in\calW\calM_0(m,N-1,N-1,\bF)$ and if we continue these
divisions, we get the factorization (\ref{pf6}).

The uniqueness of a divisor $\mbV_1(z)$ (and hence of other
factors) follows from the fact that if
$\mbV_1(z)=Q+Pz\in\calW\calM_0(m,1,1,\bF)$, where $P$ and $Q$ are
complementary projections and $\rank(P)=1$ (see Lemma 2), then
$PA_0=0$ as $\mbV_1^{-1}(z)\mbA(z)=(Pz^{-1}+Q)\mbA(z)$ does not
have the coefficient at $z^{-1}$, and since $\rank(A_0)=m-1$, such
a projection  $P$ is unique.
\end{proof}

\section{Proof of the Main Result}

We introduce the following system of $m$ conditions which plays a
key role in the proof of Theorem 1 (where this system comes from
is explained in \cite[Lemma 5]{EJL3}) . Namely, for given Laurent
polynomials
\begin{equation}
\label{zeta}
  \zeta_j(z)\in\calP_N^-[\bF],\;\;\;  j=1,2,\ldots,m-1,
 \end{equation}
 let
\begin{equation}
\label{IEEE15}
\begin{cases} \zeta_1(z)x_m(z)-\wdt{x_1}(z)\in \calP^+,\\
              \zeta_2(z)x_m(z)-\wdt{x_2}(z)\in \calP^+,\\
              \vdots\\
              \zeta_{m-1}(z)x_m(z)-\wdt{x_{m-1}}(z)\in \calP^+,\\
              \zeta_1(z)x_1(z)+\zeta_2(z)x_2(z)+\ldots+\zeta_{m-1}(z)x_{m-1}(z)
              +\wdt{x_m}(z)\in \calP^+.
         \end{cases}
\end{equation}
We say that a vector function
\begin{equation}
\label{IEEE16}
\mathbf{u}(z)=\big(u_1(z),u_2(z),\ldots,u_{m-1}(z),\wdt{u_m}(z)\big)^T,
\;\;\;u_i(z)\in \calP_N^+[\bF],\,i=1,2,\ldots,m,
\end{equation}
(we emphasize that the first $m-1$ entries of (\ref{IEEE16})
belong to  $\calP_N^+[\bF]$ and the last entry belongs to
$\wdt{\calP_N^+}[\bF]${}) is a solution of (\ref{IEEE15}) if and
only if all the conditions in (\ref{IEEE15}) are satisfied
whenever $x_i(z)=u_i(z)$, $i=1,2,\ldots,m$ (it is assumed that
$\wdt{x_m}(z)=\wdt{u_m}(z)${}).
  Observe that the set of solutions of
(\ref{IEEE15}) is a linear subspace of $\calP_N^\oplus[\bF](m):=
\underbrace{{\calP_N^+}[\bF]\tm\ldots\tm{\calP_N^+}[\bF]\tm
\wdt{\calP_N^+}[\bF]}_m$. We will see that actually this subspace
is always $m$ dimensional.

We make essential use of the following

\begin{lemma}
Let $(\ref{zeta})$  hold, and  let
\begin{equation}
\label{IEEE18a}
\mathbf{u}(z)=\big(u_1(z),u_2(z),\ldots,u_{m-1}(z),\wdt{u_m}(z)\big)^T
\in \calP_N^\oplus[\bF](m),
\end{equation}
and
\begin{equation}
\label{IEEE18b}
 \mathbf{v}(z)=\big(v_1(z),v_2(z),\ldots,v_{m-1}(z),\wdt{v_m}(z)\big)^T
 \in \calP_N^\oplus[\bF](m),
\end{equation}
be two $($possibly the same$)$ solutions of the system
$(\ref{IEEE15})$. Then
$\langle{\mathbf{u}}(z),{\mathbf{v}}(z)\rangle_m$ is the same for
every $z\in\bC\backslash\{0\}$, i.e.
\begin{equation}
\label{IEEE19}
\sum_{i=1}^{m-1}u_i(z)\wdt{v_i}(z)+\wdt{u_m}(z)v_m(z)=\const.
\end{equation}
\end{lemma}

\begin{proof}
Substituting  $x_i=v_i$ in the first $m-1$ conditions and
$x_i=u_i$ in the last condition of (\ref{IEEE15}), and then
multiplying the functions in the first $m-1$ conditions by $u_i$
and the function in the last condition by $v_m$,  we get
$$
\begin{cases}
\zeta_1v_mu_1-\wdt{v_1}u_1\in \mathcal{P}^+,\\
             \zeta_2v_mu_2-\wdt{v_2}u_2\in \mathcal{P}^+,\\
            \vdots\\
           \zeta_{m-1}v_mu_{m-1}-\wdt{v_{m-1}}u_{m-1}\in \mathcal{P}^+,\\
\zeta_1u_1v_m+\zeta_2u_2v_m+\ldots+\zeta_{m-1}u_{m-1}v_m
              +\wdt{u_m}v_m\in \mathcal{P}^+.
     \end{cases}$$

Subtracting the first $m-1$ functions from the last function in
the latter system, we get
\begin{equation}
\label{IEEE21}
\sum_{i=1}^{m-1}u_i(z)\wdt{v_i}(z)+\wdt{u_m}(z)v_m(z)\in
\mathcal{P}^+.
\end{equation}
We can interchange the roles of $u$ and $v$ in the above
discussion to get in a similar manner that
\begin{equation}
\label{IEEE22}
\sum_{i=1}^{m-1}v_i(z)\wdt{u_i}(z)+\wdt{v_m}(z)u_m(z)\in
\mathcal{P}^+.
\end{equation}
It follows from the relations (\ref{IEEE21}) and (\ref{IEEE22})
that the function in (\ref{IEEE19}) and its adjoint belong to
$\mathcal{P}^+$, which implies (\ref{IEEE19}).
\end{proof}

\begin{corollary}
If $(\ref{IEEE18a})$ and $(\ref{IEEE18b})$ are two solutions of
the system $(\ref{IEEE15})$, and
\begin{equation}
\label{1=1} \mathbf{u}(1)=\mathbf{v}(1),
\end{equation}
then
\begin{equation}
\label{z=z} \mathbf{u}(z)=\mathbf{v}(z) \;\text{ for each
}z\in\bC\backslash\{0\}.
\end{equation}
In particular, if $\mathbf{u}(1)=\mbo\in\bF^m$, then
$\mbu(z)=\mbo$ for each $z\in\bC\backslash\{0\}$ since the trivial
vector function $\mbv(z)=\mbo$ is always a solution of the system
$(\ref{IEEE15})$.
\end{corollary}

\begin{proof}
 Since $\mathbf{u}(z)-\mathbf{v}(z)$ is also a solution of
(\ref{IEEE15}), it follows from the lemma that
$\|{\mathbf{u}}(z)-{\mathbf{v}}(z)\|_{m}=
\|{\mathbf{u}}(1)-{\mathbf{v}}(1)\|_{m}=0$ for each
$z\in\bC\backslash\{0\}$. Hence (\ref{z=z}) holds.
\end{proof}

\begin{remark}
One can see that Corollary $1$ remains valid if we take any fixed
point $z_0\not=0$ in the role of $1$ in the equation
$(\ref{1=1})$.
\end{remark}

\begin{lemma}
Let $(\ref{zeta})$  hold, and  let vector functions
\begin{equation}
\label{IEEE38}
\mathbf{v}_j(z)=\big(v_{1j}(z),v_{2j}(z),\ldots,\wdt{v_{mj}}(z)\big)^T\in
\calP_N^\oplus[\bF](m),\;\;\;j=1,2,\ldots,m,
\end{equation}
be any $m$ solutions of the system $(\ref{IEEE15})$. Then the
determinant of the $m\tm m$ Laurent polynomial matrix
\begin{equation}
\label{IEEE381} V(z)=({\mathbf{v}_1}(z),{\mathbf{v}_2}(z),\ldots
{\mathbf{v}_m}(z))
\end{equation}
is constant
\begin{equation}
\label{IEEE41} \det V(z)=\const., \;\;\;z\in\bC\backslash\{0\}.
\end{equation}
\end{lemma}
\begin{proof}
Since the vector polynomials (\ref{IEEE38}) satisfy the last
condition of the system  (\ref{IEEE15}), we have
\begin{equation}
\label{325} F(z)V(z)\in\calP^+(m\tm m),
\end{equation}
where $F(z)$ is defined by (\ref{F93}). Consequently, as $\det
F(z)=1$,
 \begin{equation}
 \label{IEEE43}
\det V(t)\in \mathcal{P}^+.
\end{equation}
Since the vector polynomials (\ref{IEEE38}) satisfy the first
$m-1$ conditions of the system  (\ref{IEEE15}), we have
$$
\phi_{ij}(z):=\zeta_i(z)v_{mj}(z)-\wdt{v_{ij}}(z)\in\calP^+,
$$
$1\leq j\leq m$, $1\leq i< m$. Direct computations show that
$$
\wdt{V}(z)F^{-1}(z)=\Phi(z)\in\calP^+(m\tm m),
$$
where $F^{-1}(z)$ is obtained from $F(z)$ by replacing each
$\zeta_i$ by $-\zeta_i$ and
$$
\Phi(z)= \left(\begin{matrix}
           -{\phi_{11}(z)}&-{\phi_{21}(z)}&\cdots&-{\phi_{m-1,1}(z)}& {v_{m1}(z)}\\
           -{\phi_{12}(z)}&-{\phi_{22}(z)}&\cdots&-{\phi_{m-1,2}(z)}& {v_{m2}(z)}\\
                       \vdots&    \vdots&\vdots&\vdots&\vdots\\
  -{\phi_{m1}(z)}&-{\phi_{m1}(z)}&\cdots&-{\phi_{m-1,m}(z)}&{v_{m,m}(z)}\\
           \end{matrix}\right).
$$
Consequently $\det \wdt{V}(z)=\det\Phi(z)\in\calP^+$ and together
with (\ref{IEEE43}) this gives (\ref{IEEE41}).
\end{proof}

\begin{remark} Observe that the relation $(\ref{325})$ holds under
the hypothesis of the lemma.
\end{remark}

Now we construct explicitly $m$ independent solutions
(\ref{IEEE38}) of the system (\ref{IEEE15}). We seek for a
nontrivial  solution
\begin{equation}
\label{v2_46}
\mathbf{x}(z)=\big(x_1(z),x_2(z),\ldots,\wdt{x_m}(z)\big)^T\in
\calP_N^\oplus[\bF](m)
\end{equation}
of (\ref{IEEE15}), where
\begin{equation}
\label{v2_47}
 x_i(z)=\sum_{k=0}^N \al_{ik} z^k,\;\;\;i=1,2,\ldots,
m,
\end{equation}
and the coefficients $a_{ik}$ are to be determined.

Equating the coefficients of $z^{-k}$, $k=0,1,2,\ldots,N$ of the
Laurent polynomials in the system (\ref{IEEE15}) to 0, except the
$0$th coefficient of the $j$th function which we equate to 1, we get
the following system of algebraic equations in the block matrix form
which we denote by $\mathbb{S}_j$:
\begin{equation}
\label{v2_48}
 \mathbb{S}_j:= \begin{cases}\Tht_1 X_m-\overline{X_1}={\bf 0}, \\
    \Tht_2 X_m-\overline{X_2}={\bf 0}, \\
    \vdots    \\
    \Tht_j X_m-\overline{X_j}={\bf 1}, \\
    \vdots    \\
    \Tht_{m-1} X_m-\overline{X_{m-1}}={\bf 0}, \\
    \Tht_1 X_1+\Tht_2 X_2+\ldots+\Tht_{m-1}
    X_{m-1}+\overline{X_m}={\bf 0}\;, \end{cases}
\end{equation}
where $\Tht_i$ is defined from the equations (\ref{tht1}) and
(\ref{zeta1}), ${\bf 0}=(0,0,\ldots,0)^T\in \mathbb{F}^{N+1}$,
${\bf 1}=(1,0,0,\ldots,0)^T\in \mathbb{F}^{N+1}$, and the column
vectors
\begin{equation}
\label{v2_50}
X_i=(\al_{i0},\al_{i1},\ldots,\al_{iN})^T,\;\;i=1,2,\ldots,m,
\end{equation}
(see (\ref{v2_47}){}) are unknown.

\begin{remark}
We emphasize  that if $(X_1,X_2,\ldots,X_m)$ defined by
$(\ref{v2_50})$ is the solution of the system $(\ref{v2_48})$,
then the vector function $(\ref{v2_46})$ defined by
$(\ref{v2_47})$ will be the solution of the system
$(\ref{IEEE15})$.
\end{remark}

It is easy to show that the system (\ref{v2_48}) $\mathbb{S}_j$
has a solution for each $j=1,2,\ldots,m$. Indeed, determining
$X_i$, $i=1,2,\ldots,m-1$, from the first $m-1$ equations in
(\ref{v2_48}),
\begin{equation}
\label{v2_51}
X_i=\overline{\Tht_i}\cdot\overline{X_m}-\delta_{ij}{\bf
1},\;\;i=1,2,\ldots,m-1,
\end{equation}
 and then substituting  into the last equation of (\ref{v2_48}), we get
$$
\Tht_1\cdot\overline{\Tht_1}\cdot\overline{X_m}+\Tht_2\cdot\overline{\Tht_2}\cdot\overline{X_m}
+\ldots+\Tht_{m-1}\cdot\overline{\Tht_{m-1}}\cdot\overline{X_m}+\overline{X_m}=
\Tht_j\cdot{\bf 1}
$$
(we assume that $\Tht_m=I_{N+1}$, i.e. the right-hand  side is
equal to ${\bf 1}$ when $j=m$) or, equivalently,
\begin{equation}
\label{v2_52}
(\Theta_1\cdot\overline{\Theta_1}+\Theta_2\cdot\overline{\Theta_2}
+\ldots+\Theta_{m-1}\cdot\overline{\Theta_{m-1}}+I_{N+1})\cdot\overline{X_m}
=\Tht_j\cdot{\bf 1}.
\end{equation}
This is the same system as (\ref{517}) (see \ref{Dt1}). Since each
$\Tht_i$ is symmetric, $\Tht_i=\Tht_i^T$, we have that
$\Theta_i\overline{\Theta_i}=\Theta_i\Theta_i^*$,
$i=1,2,\ldots,m-1$, are non-negative definite and the coefficient
matrix (\ref{Dt1}) of the system (\ref{v2_52}) (which is the same
for each $j=1,2,\ldots,m$) is positive definite (with all
eigenvalues larger than or equal to 1). Consequently, $\Delta$ is
nonsingular, $\det\Delta\geq1$, and the system (\ref{v2_52}) has a
unique solution for each $j$.

Finding the  vector $\overline{X_m}\in\bF^{N+1}$ from
(\ref{v2_52}) and then determining $X_1,X_2,\ldots,X_{m-1}$ from
(\ref{v2_51}), we get the unique solution of $\mathbb{S}_j$. To
indicate its dependence on $j$, we denote the solution of
$\mathbb{S}_j$ by $(X_1^j,X_2^j,\ldots,X_{m-1}^j,X_m^j)$,
\begin{equation}
\label{v2_54} X_i^j:=(\al_{i0}^j,\al_{i1}^j,\ldots, \al_{iN}^j)^T,
\;\;\;i=1,2,\ldots,m,
\end{equation}
so that the vector functions (\ref{IEEE38}), where
$$
v_{ij}(z)=\sum_{k=0}^N \al_{ik}^j z^k, \;\;\;1\leq i,j\leq m,
$$
are $m$ solutions of the system (\ref{IEEE15}) (see Remark 5).
These vector functions $\mbv_1(z)$, $\mbv_2(z)$, $\ldots$,
$\mbv_m(z)$ are linearly independent since the linear
transformation ${\bf L}:\calP_N^\oplus[\bF](m)\to\bF^m$ which maps
$(x_1(z),x_2(z), \ldots,\wdt{x_m}(z))$ into the $0$th coefficients
of the functions (or their adjoints) standing on the left-hand
side of the system (\ref{IEEE15}), viz. into
$\big(c_0\{\wdt{\zeta_1}(z)\wdt{x_m}(z)-{x_1(z)}\}$,
$c_0\{\wdt{\zeta_2}(z)\wdt{x_m}(z)-{x_2(z)}\}$, $\ldots$,
$c_0\{\zeta_1(z)x_1(z)+\zeta_2(z)x_2(z)+\ldots+\wdt{x_m}(z)\}\big)
$, transforms $m$ vector functions $\mbv_1(z)$, $\mbv_2(z)$,
$\ldots$, $\mbv_m(z)$  into linearly independent standard bases of
$\mathbb{F}^m$. Namely, ${\bf L}\big(\mbv_i(z)\big)
=(\delta_{i1},\delta_{i2},\ldots,\delta_{im})$, $i=1,2,\ldots,m$,
because of (\ref{v2_48}). Consequently $V(1)$ is nonsingular (see
(\ref{IEEE381}){}) since if ${\bf w}\cdot V(1)={\bf 0}\in\bF^m$
for some ${\bf 0}\not={\bf w}=(w_1,w_2,\ldots,w_m)
\in\mathbb{F}^m$, then ${\bf w}\cdot
V(z)=\sum_{k=1}^mw_k\mbv_k(z)={\mbo}$ by virtue of Corollary 1 of
Lemma 3, which contradicts  the independence of $\mbv_1(z)$,
$\mbv_2(z)$,$\ldots$,$\mbv_m(z)$.

Let $U(z)$ be defined by (\ref{V-1}). Then it has the form
(\ref{U6}) and its column vectors $\mbu_1(z)$,
$\mbu_2(z)$,$\ldots$,$\mbu_m(z)$,
\begin{equation}
\label{VU} U(z)=({\mathbf{u}_1}(z),{\mathbf{u}_2}(z),\ldots,
{\mathbf{u}_m}(z)),
\end{equation}
are $m$ solutions of (\ref{IEEE15}) satisfying
\begin{equation}
\label{uiei} \mbu_i(1)={\bf e}_i,\;\;\;i=1,2,\ldots,m,
\end{equation}
 (since (\ref{UIm4}) holds because
of (\ref{V-1}){}). Therefore, the matrix function $U(z)$ is
paraunitary since $\wdt{U}(z)U(z)$ is a constant matrix, by virtue
of Lemma 3, and this constant is $I_m$ since (\ref{UIm4}) holds.
By virtue of Lemma 4, the determinant of $U(z)$ is also a constant
which is equal to 1 since $\det u(1)=1$. Hence the relations
(\ref{UPp1})--(\ref{UIm4}) are valid for (\ref{VU}). Observe also
that (\ref{Pr95}) holds because of (\ref{325}) (see Remark 4) and
(\ref{V-1}). Hence it remains to show only the relation
(\ref{Ulr5}) in order to complete the proof of the first part of
Theorem 1. Meanwhile we are ready to formulate and prove a
corollary of the above discussion which will be used in the next
section.

\begin{corollary}
Let $(\ref{zeta})$ hold, let $F(z)$ be the matrix function defined
by $(\ref{F93})$, and let $(\ref{VU})$ be the corresponding
$($according to Theorem $1${}$)$
$U(z)\in\mathcal{P}\mathcal{U}_1(m,N,\bF)$.

If ${\bf b}=(b_1,b_2,\ldots,b_m)\in\bF^m$, then $\mbu_{\bf b}(z)$
is a solution of the system $(\ref{IEEE15})$  satisfying
$$
\mbu_{\bf b}(1)={\bf b},
$$
if and only if
$$ \mbu_{\bf
b}(z)=\sum_{i=1}^mb_i\,\mbu_i(z),\;\;\;z\in\bC\backslash\{0\}.
$$
\end{corollary}

\begin{proof}
Since $\mathbf{u}_1(z),{\mathbf{u}_2}(z),\ldots,
{\mathbf{u}_m}(z)$ (see (\ref{VU}){}) are  solutions of the system
(\ref{IEEE15}) (recall that the set of its solutions is a linear
subspace) and (\ref{uiei}) holds, the first part of the corollary
is clear. The second part follows from the first part and from
Corollary 1.
\end{proof}

In order to prove (\ref{Ulr5}), let
\begin{equation}
\label{Psi} \Psi(z)=\big(\psi_{ij}(z)\big)_{i,j=1}^m:=F(z)U(z),
\end{equation}
 which, as  has already been mentioned,
 belongs to $\calP^+(m\tm m)$ (see (\ref{Pr95}){}). Since the  determinants of $F(z)$ and $U(z)$ are 1
for each $z\in\bC\backslash\{0\}$, we can conclude that
\begin{equation}
\label{dtpsi} \det\Psi(z)=1\;\;\text{ for each } z\in\bC.
\end{equation}
Because of the structure of the matrix (\ref{F93}), we have
$$
\psi_{ij}(z)=u_{ij}(z),\;\;\;1\leq i<m,\;1\leq j\leq m.
$$
Since $\wdt{U}(z)=U^{-1}(z)$ and (\ref{Udet3}) holds, we have
$$
u_{mj}(z)=\cof u_{mj}(z)=\cof \psi_{mj}(z),\;\;1\leq j\leq m.
$$
Thus, by virtue of (\ref{dtpsi}),
$$
1=\det\Psi(z)=\sum_{j=1}^m\psi_{mj}(z)\cof\psi_{mj}(z)=
\sum_{j=1}^m\psi_{mj}(z)u_{mj}(z)
$$
for each $z\in\bC$ including $0$. Consequently (\ref{Ulr5}) holds.

Thus, the constructed matrix function $U(z)$ has all the desired
properties, i.e. (\ref{U94}) and (\ref{Pr95}) hold. Note that a
brief way of construction of the entries of $U(z)$ is described by
the formulas (\ref{zeta1})--(\ref{V-1}).

The uniqueness of (\ref{U94}) follows from the uniqueness of
spectral factorization (see Section 2) since $F(z)U(z)$ is the
spectral factor of $F(z)\wdt{F}(z)$.

Let us now show the converse part of Theorem 1. If we have a
matrix function (\ref{U94}), then $u_{mj}(0)\not=0$ for some $j$
since (\ref{Ulr5}) holds. Let us determine functions $\zeta_i$ by
the formula (\ref{zeta24}), construct the matrix function
(\ref{F93}), and show that (\ref{Pr95}) is valid. For this we need
only to check that the entries of the last row of $\Psi(z)$ (see
(\ref{Psi}){}) belong to $\calP^+[\bF]$. For $1\leq n\leq m$, we
have
\begin{gather}
\calP\ni\sum_{i=1}^{m-1}\zeta_i(z)u_{in}(z)+\wdt{u_{mn}}(z)=
\sum_{i=1}^{m-1} \left(\frac{\wdt{u}_{ij}(z)}{u_{mj}(z)}
-\left[\frac{\wdt{u}_{ij}(z)}{u_{mj}(z)}\right]^+\right)
u_{in}(z)+\wdt{u_{mn}}(z)=\notag\\\label{shua}
\frac{1}{u_{mj}(z)}\left(\sum_{i=1}^{m-1}
\wdt{u}_{ij}(z)u_{in}(z)+u_{mj}(z)\wdt{u_{mn}}(z)\right)-
\sum_{i=1}^{m-1} \left[\frac{\wdt{u}_{ij}(z)}{u_{mj}(z)}\right]^+
u_{in}(z)=\\
\frac{\dt_{nj}}{u_{mj}(z)}- \sum_{i=1}^{m-1}
\left[\frac{\wdt{u}_{ij}(z)}{u_{mj}(z)}\right]^+ u_{in}(z),\notag
\end{gather}
which is analytic in a neighborhood of $0$ and yields that this
function belongs to $\calP^+$. Thus (\ref{Pr95}) holds.

Let us now show the uniqueness of the desired $F(z)$.

Recall that  $F^{-1}(z)$ can be  obtained  by replacing each
$\zeta_i$ by $-\zeta_i$ in (\ref{F93}). Hence
$F^{-1}(z)\in\wdt{\calP_N^+}[\bF](m\tm m)$. Since
$\Psi(z)\in\calP^+(m\tm m)$ (see (\ref{Psi}) and (\ref{Pr95}){})
and (\ref{dtpsi}) holds, we have $\Psi^{-1}(z)\in\calP^+(m\tm m)$.

If $F_1(z)$ has the same form (\ref{F93}) with the last row
$[\zeta_1'(z),\zeta_2'(z),\ldots,\zeta_{m-1}'(z),1]$ and
$\Psi_1(z):=F_1(z)U(z)\in\calP^+(m\tm m)$, then
$U(z)=F^{-1}(z)\Psi(z)=F_1^{-1}(z)\Psi_1(z)$, which yields
$$
\wdt{\calP^+}(m\tm m)\ni
F_1(z)F^{-1}(z)=\Psi_1(z)\Psi^{-1}(z)\in\calP^+(m\tm m).
$$
Thus $F_1(z)F^{-1}(z)$ is a constant matrix while, on the other hand,
this product has the form (\ref{F93}) with the last row
$[\zeta_1'-\zeta_1,\zeta_2'-\zeta_2,\ldots,\zeta_{m-1}'-\zeta_{m-1},1]$,
 which implies that $\zeta'_i=\zeta_i$, $i=1,2,\ldots,m-1$.

 \section{Wavelet Matrices with Given First Rows}

 If we know the first row of the wavelet matrix (\ref{pf6}), then one
 can uniquely determine each $P_i$ and $Q_i$ in (\ref{pf6}) and
 then recover $\mbA(z)$ by the formula (\ref{pf6}), i.e. find the
 remaining $m-1$ rows of $\mbA(z)$. The algorithm for this
 procedure is well known (see e.g. \cite[Th. 4.4.17]{RW}). We
 describe a new method of reconstruction of $\mbA(z)$ based on the
 proposed parametrization of compact wavelet matrices. This
 description is included in a constructive proof of the following

 \begin{theorem}
Let ${\rm v}=(v_1,v_2,\ldots,v_m)\in\bF^m$ be a vector of the unit
norm and $V\in\bF^{m\tm m}$ be a unitary matrix with the first row
${\rm v}$. If
$\mbA_1(z)=\big({a}_{11}(z),{a}_{12}(z),\ldots,{a}_{1m}(z)\big)\in\calP^+_N(1\tm
m)$ is a row function of order $N$
$\big(${}$\sum_{j=1}^m|\al_{Nj}|>0$, where
${a}_{1j}(z)=\sum_{k=0}^N\al_{kj}z^k$, $j=1,2,\ldots,m${}$\big)$
such that $\mbA_1(z)\wdt{\mbA_1}(z)=I_m$, $z\in\bC\backslash\{0\}$
and
\begin{equation}
\label{A1v} {\mathbf A}_1(1)={\rm v},
\end{equation}
then there exists a unique wavelet matrix of rank $m$ and of order
and degree $N$, $\mbA(z)\in\calW\calM(m,N,N,\bF)$, with
$\mbA(1)=V$ which has the first row $\mbA_1(z)$.

\begin{remark}
For a given ${\rm v}\in\bF^m$ with $\|{\rm v}\|_m=1$, it is in
general hard to decide whether there exists a unitary matrix
$V\in\bF^{m\tm m}$ with the first row $v$ $($because of the
arbitrariness of the field $\bF$ which not always contains
$\sqrt{a}$ for a positive $a\in\bF${}$)$. However the theorem can
be always used for ${\rm v}={\bf e}_1$ and $V=I_m$.
\end{remark}

\begin{proof}
For notational convenience, we will find $\mbA^T(z)$ instead of
$\mbA(z)$. Thus we assume that the ${\rm
v}=(v_1,v_2,\ldots,v_m)^T$ is a column vector and
$\mbA_1(z)=\big({a}_{11}(z),\ldots,{a}_{1m}(z)\big)^T\in\calP^+_N(m\tm
1)$, which satisfies $\wdt{\mbA_1}(z){\mbA_1}(z)=I_m$ and
(\ref{A1v}). Furthermore, without loss of generality (we
interchange the rows if necessary), we assume that
$\al_{Nm}\not=0$ where we recall
${a}_{1m}(z)=\sum_{k=0}^N\al_{km}z^k$.

Consider now
$\mbU_1(z)=\big(u_{11}(z),u_{21}(z),\ldots,\wdt{u_{m1}}(z)\big)^T
\in\calP^\oplus_N(m)$, where $u_{j1}(z)=a_{1j}(z)$,
$j=1,2,\ldots,m-1$, and $\wdt{u_{m1}}(z)=z^{-N}a_{1m}(z)$. Then
$u_{m1}(0)\not=0$, and we can define $\zeta_i$,
$i=1,2,\ldots,m-1$, by a formula like (\ref{zeta24})
$$
\zeta_i(z)=\left[\frac{\wdt{u}_{i1}(z)}{u_{m1}(z)}\right]^-=
\left[\frac{\wdt{u}_{i1}(z)}{u_{m1}(z)}\right]_N^-,
\;\;\;\;i=1,2,\ldots,m-1,
$$
and construct the matrix function $F(z)$ defined by (\ref{F93}).
Consider now the corresponding (according to Theorem 1)
$U(z)\in\mathcal{P}\mathcal{U}_1(m,N,\bF)$ whose columns
$\mbu_1(z)$, $\mbu_2(z)$, $\ldots$,$\mbu_m(z)$, as we remember,
are solutions of the system (\ref{IEEE15}) and also satisfy
(\ref{uiei}). We can make sure by direct computations like
(\ref{shua}) that $U_1(z)$ satisfies the last condition of
(\ref{IEEE15}), and taking into account the relations
$$
\zeta_i(z)u_{m1}(z)-\wdt{u_{i1}}(z)=\left(\frac{\wdt{u}_{i1}(z)}{u_{m1}(z)}
-\left[\frac{\wdt{u}_{i1}(z)}{u_{m1}(z)}\right]^+\right)
u_{m1}(z)-\wdt{u_{i1}}(z)=-\left[\frac{\wdt{u}_{i1}}{u_{m1}}\right]^+
u_{m1}\in\calP^+,
$$
$i=1,2,\ldots,m-1$, we see that $U_1(z)$ is a solution of the
system (\ref{IEEE15}). Consequently, according to Corollary 2,
$$
U_1(z)=U(z)\cdot {\rm v}=\sum_{i=1}^mv_i\,\mbu_i(z),
$$
(as this equation holds for $z=1$) and it is the first column of
$U(z)\cdot V$. This matrix function is paraunitary as well (since
$V$ is unitary) and satisfies the condition $U(1)V=V$. Hence
$$\mbA(z)=\diag[1,1,\ldots,1,z^N]U(z)V$$
 will be the desired
matrix.

The uniqueness of $\mbA(z)$ is also valid since $F(z)$ described in
the above proof  is the same as the matrix function corresponding to
$\mbA(z)\cdot V^{-1}\in\calW\calM_1(m,N,N,\bF)$
by the to one-to-one map (\ref{W1U1}) and Theorem 1, and such
$F(z)$ is unique.
\end{proof}

 \end{theorem}

\section{Other Possible Applications}

In the end,  we would like to discuss briefly some possible
applications of the proposed method of compact wavelet matrices
construction. It should be emphasized that the intent of this
section is largely motivational and we  consider three separate
topics.

{\bf A. Rational approximations to compact wavelet matrices.} Most
of compact wavelet matrices used in practice, for example
Daubechies  wavelet matrices, obey certain additional
restrictions. The proofs of the existence of such matrices and the
ways of their construction are highly non-linear and thus the obtained
coefficients are in general irrational.  In actual calculations on a
digital computer, these coefficients should be quantized and hence
approximated by rational numbers. It may happen during this
approximation that the basic property of the wavelet matrices
(\ref{Qc2}) will not be preserved exactly. Using the proposed
parametrization and taking $\bQ$ in the role of $\bF$, we can
provide an approximation of any compact wavelet matrix  $\calA$ by
$\calA'$ with rational coefficients preserving exactly the shifted
orthogonality condition (\ref{Qc2}). Indeed, it is just
sufficient to note  in the proof of Theorem 1 that the maps
described in  (\ref{PAR}) are continuous. Hence, if we find the point in $\bR^{(m-1)N}$
corresponding to $\calA$, approximate it
by rational coordinates  and go back into the space of wavelet
matrices, then we will get the desired $\calA'$ (in a similar
manner, Vaidyanathan \cite{Vai} used the parametrization
(\ref{par1}) for the same purposes). Such various  approximations
for Daubechies wavelet matrices of different genus are explicitly
constructed in \cite{EGL}.

{\bf B. Cryptography.} As  said in Remark 1, the factorization
(\ref{pf6}) much resembles  the factorization for ordinary
polynomials into linear terms. On the other hand, an efficient way
of construction of paraunitary matrix polynomials associated to
the given points in $\bF^{(m-1)N}$  expressed by the diagram
(\ref{PAR}), which helps to handle such matrix polynomials easily,
can be compared to natural parametrization of ordinary polynomials
by to their coefficients. An application of polynomial
factorization theory in the cryptography is widely known. Having
the above similarities (and the advantages mentioned in Remark 1),
one can also expect certain applications of the developed theory
in the cryptography.

{\bf C. Selection of best wavelet matrices.} Which wavelet matrix
is most suitable  to apply in a given practical situation
represents frequently a certain optimization problem (which is
sometimes very hard to solve) or should be obtained empirically by
computer simulations. Having a quick access to the complete bank
of compact wavelet matrices due to the parametrization (\ref{PAR})
gives an opportunity to choose the best possible wavelet matrix by
nearly complete screening.

 \vskip+0.2cm

\ Authors' Addresses: \vskip+0.2cm

 L. Ephremidze, E. Lagvilava

\noindent  A. Razmadze mathematical Institute

\noindent I. Javakhishvili State University

\noindent 2, University Street, Tbilisi 0143, Georgia

\noindent E-mails: {\em lephremi@umd.edu; edem@rmi.ge}

\end{document}